\documentclass{article}
\usepackage{amsmath,amsgen,amstext,amsbsy,amsopn,amsfonts}
\newtheorem{Def}{Definition}[section]
\newtheorem{Thm}[Def]{Theorem}
\newtheorem{Prp}[Def]{Proposition}
\newtheorem{Lem}[Def]{Lemma}
\newtheorem{Cor}[Def]{Corollary}

\newcommand{\g}{\mathfrak{g}}

\newcommand{\half}{\frac{1}{2}}

\begin{document}

\centerline{\LARGE Splints of classical root systems}

\vskip .25 in
\begin{flushleft}
David A. Richter \\
Department of Mathematics, MS 5248 \\ 
Western Michigan University \\ 
Kalamazoo MI 49008-5248 \\
\texttt{david.richter@wmich.edu} \\
\end{flushleft}

{\bf Abstract.}
This article introduces a new term ``splint'' and classifies the splints of the
classical root systems.  The motivation comes from representation theory
of semisimple Lie algebras.  In a few instances, splints play a role in determining branching
rules of a module over a complex semisimple Lie algebra when restricted to a subalgebra.  In these particular
cases, the set of submodules with respect to the subalgebra themselves may be regarded as the character
of another Lie algebra.

\noindent
Keywords:  branching rule; root system; splint

\noindent
AMS Mathematics Subject Classification (2000):  17B05, 17B10.

\section{Introduction}

This paper is concerned with a branching rule for a pair $\g_1\subset\g$ of semisimple Lie algebras
in which the highest weights of a representation $V$ with respect to $\g_1$ themselves may serve as the weights
of another semisimple Lie algebra.  The particular phenomenon has been mentioned in the literature
only a few times, and generally only parenthetically; see \cite{Hec} and \cite{KQ}.

It is perhaps easiest to see the phenomenon the case where $\g$ is the
14-dimensional exceptional simple Lie algebra $\g_2$ and the subalgebra is $\mathfrak{sl}_3$.
If $V$ is any finite-dimensional irreducible module for $\g_2$, then the branching into
irreducible representations of $\mathfrak{sl}_3$ has the pattern of a weight diagram for $\mathfrak{sl}_3$.
As \cite{Hec} explains, this follows from the fact that the difference between the root systems
of $\g_2$ and $\mathfrak{sl}_3$ is a copy of the root system for $\mathfrak{sl}_3$.

A less trivial example occurs when
$\g=\mathfrak{b}_2\cong\mathfrak{so}_5$.  Recall that the set
$$\Delta^+=\{e_1,e_2,e_1-e_2,e_1+e_2\}$$
serves as a system of positive roots for this Lie algebra, where $\{e_1,e_2\}$ is an orthonormal
basis for the dual of a prescribed Cartan subalgebra.
Choose $\g_1=\mathfrak{a}_1\cong\mathfrak{sl}_2$ as
the Lie subalgebra corresponding to the subset
$$\Delta_1^+=\{e_1-e_2\}.$$
Next, let $V$ be a finite-dimensional irreducible module over $\g$.  Then $V$ is
a module over $\g_1$ and one may write
$$V=\bigoplus_\mu V_\mu,$$
where each $V_\mu$ is an irreducible $\g_1$-module with highest weight $\mu$.  
Then, as the reader may check, the weights $\mu$ appearing in this summand 
have the pattern of a weight diagram of a representation of $\mathfrak{sl}_3$.  This occurs despite
the fact that $\mathfrak{a}_2$ is not a subalgebra of $\mathfrak{b}_2$.

These examples motivate the introduction of the term ``splint''.  Roughly, a splint of a root system
$\Delta$
is a partition $\Delta_1\cup\Delta_2$ into two subsets, each of which have the additive, but not necessarily
metrical, properties of a root system.
The aim in this study is to classify all splints of root systems for complex semisimple Lie algebras.
This problem is slightly more general
than what is needed with regards to branching rules.  In the example above, one of the ``components''
$\Delta_1$ corresponds to a Lie subalgebra of $\g$, but one will not find this constraint in the definitions;
there are few enough cases to sort through if the reader wishes to find these.

One of the ingredients of the term ``splint'' is the idea of a root system embedding.
The definition given here is novel, although \cite{BP} introduces a similar definition.
The reference \cite{LM} is given because the techniques used there seem the closest to the ideas expressed
in this introduction.

\section{Definitions}

The purpose of this section is is to define
what it means for a root system to splinter, but one first needs some preliminary terminology.
First of all,
the term ``simple root system'' refers to one appearing in Cartan's list
$$A_r\ (r\geq 1),\ B_r\ (r\geq 2),\ C_r\ (r\geq 3),\ D_r\ (r\geq 4),\ E_6,\ E_7,\ E_8,\ F_4,\ G_2,$$
corresponding to four classical sequences of and five exceptional complex simple Lie algebras, \cite{Hum,Jac}.

If $\g$ is complex semisimple, then its root system is a disjoint union
of various copies of these, each corresponding to a summand of $\g$ into
simple components.  One may abstract this as follows:  Suppose $\Delta_1$ and $\Delta_2$ are root systems.
Then one assumes already that $\Delta_1$ and $\Delta_2$ are subsets of vector spaces $V_1$ and $V_2$ respectively.
The ``direct sum'' of the root systems is then
$$\Delta_1+\Delta_2=\{(\alpha,0):\alpha\in\Delta_1\}\cup\{(0,\alpha):\alpha\in\Delta_2\}\subset V_1\oplus V_2.$$
The direct sum operation has the following features:  First, if $\Delta=\Delta_1+\Delta_2$ with $\alpha_1\in\Delta_1$
and $\alpha_2\in\Delta_2$, then $\alpha_1+\alpha_2\notin\Delta$.  Second, if there exist roots $\alpha_1,\alpha_2\in\Delta$
such that $\alpha_1+\alpha_2$ is also a root, then the root system cannot be written as a direct sum
$\Delta=\Delta_1+\Delta_2$ where $\alpha_1\in\Delta_1$ and $\alpha_2\in\Delta_2$.
Third, if one may write $\Delta=\Delta_1\cup\Delta_2$ as a disjoint union where for all $\alpha_1\in\Delta_1$
and all $\alpha_2\in\Delta_2$, the vector $\alpha_1+\alpha_2$ is not a root, then $\Delta$ is 
the direct sum of $\Delta_1$ and $\Delta_2$.

\begin{Def}
Suppose $\Delta_0$ and $\Delta$ are root systems.  Then $\iota$ is an ``embedding''
if (a) it is a 1-1 function from $\Delta_0$ to $\Delta$, and (b)
$$\iota(\gamma)=\iota(\alpha)+\iota(\beta)$$
for all $\alpha,\beta,\gamma\in\Delta_0$ such that $\gamma=\alpha+\beta$.
\end{Def}

It is convenient to write ``$\Delta_0\hookrightarrow\Delta$'' in order to imply that there is
an embedding from $\Delta_0$ to $\Delta$ and to say, under these circumstances,
that $\Delta_0$ ``embeds'' in $\Delta$.  Given this definition, one may summarize
some of the preceeding observations:

\begin{Prp}
Suppose $\Delta$ is a root system.  Then at least one simple component of $\Delta$ is not isomorphic to $A_1$ 
if and only if there is an embedding $A_2\hookrightarrow\Delta$.
\end{Prp}

The following tells what it means for a root system to splinter.

\begin{Def}
A root system $\Delta$ ``splinters'' as $(\Delta_1,\Delta_2)$ if there are two embeddings
$\iota_1:\Delta_1\hookrightarrow\Delta$ and $\iota_2:\Delta_2\hookrightarrow\Delta$
where (a) $\Delta$ is the disjoint union of the images of $\iota_1$ and $\iota_2$ and
(b) neither the rank of $\Delta_1$ nor the rank of $\Delta_2$ exceeds the rank of $\Delta$
\end{Def}

It is equivalent to say that $(\Delta_1,\Delta_2)$ is a ``splint'' of $\Delta$.  
Each component $\Delta_1$ and $\Delta_2$ is a ``stem'' of the splint $(\Delta_1,\Delta_2)$.

Thus far, these definitions do not mention inner products.  As explained earlier, since
the factorization of the Weyl denominator does not require an inner product,
the focus here is on the additive properties of root systems.  Nevertheless,
as defined in many places, for example as in \cite{Hum}, a root system is a particular type of
subset of a vector space equipped with a Euclidean inner product.  In this regard, it is useful
to distinguish between ``metric root systems,'' which are endowed with an inner product and mere
``root systems'', described above.

It is also useful to introduce some terminology relating metric root systems with general root systems.
Suppose $\iota:\Delta_0\hookrightarrow\Delta$
is an embedding, and suppose that the inner products associated to $\Delta_0$ and $\Delta$, respectively, are
$\left<\ ,\ \right>_0$ and $\left<\ ,\ \right>$.   Call the embedding $\iota$ ``metric'' if 
there is a positive scalar $\lambda$ such that
$$\left<\alpha,\beta\right>_0=\lambda\left<\iota(\alpha),\iota(\beta)\right>$$
for all $\alpha,\beta\in\Delta_0$ and ``non-metric'' otherwise.  An imbedding 
$\iota:\Delta_0\hookrightarrow\Delta$ 
may not be metric even if the restriction of $\iota$ to every simple component is metric; call such
a non-metric embedding ``semimetric''.

The structure of a  root system is characterized by the additive properties of its corresponding
system of positive roots.  Thus, if $\Delta$ is the set of all roots, positive and negative, 
and $\Delta_1^+$ and $\Delta_2^+$ are two
choices of positive roots, then there is a bijection $f:\Delta_1^+\rightarrow\Delta_2^+$ having the property that
$f(\alpha+\beta)$ is a root if and only if $\alpha+\beta$ is a root. 
Hence, it is necessary only to consider {\it positive} root systems.  This should be considered as the 
0th definition of this section:  A 
``simple root system'' is any system of positive roots for one of the root systems listed above, and, more 
generally, a ``root system'' is any direct sum of these.

In the sequel, it is necessary to be familiar with some presentations of simple root systems.
Set
$$A_r=\{e_i-e_j:0\leq i<j\leq r\},\ D_r=\{e_i\pm e_j:1\leq i<j\leq r\},$$
$$B_r=D_r\cup\{e_i:1\leq i\leq r\},\hbox{ and }C_r=D_r\cup\{2e_i:1\leq i\leq r\},$$
where $\{e_i:i\in I\}$ is assumed to be orthonormal for any set $I$.
Notice that every root of $A_r$ and every root of $D_r$ has norm 2.  
None of the well-known presentations of $E_6$, $E_7$, and $E_8$ are required here, but one
does need the convention that every element of any of these root system has norm 2.
Finally,
$$F_4=\{2e_i:1\leq i\leq 4\}\cup\{e_1\pm e_2\pm e_3 \pm e_4\} \cup \{e_i\pm e_j:1\leq i<j\leq 4\}$$
and
$$G_2=\{e_1-e_2,e_2-e_3,e_1-e_3,2e_1-e_2-e_3,e_1+e_2-2e_3,e_1-2e_2+e_3\}.$$

\section{Results}

The main result is the following:

\begin{Thm}
The table of splints of simple root systems given in the appendix is complete.
\end{Thm}

In order to establish this, one first needs some moderately general results.  
First, it should be clear that the most difficult part of the problem
involves splintering the simple root systems.  Indeed, if a root system is a direct sum 
$\Delta=\Delta_1+\Delta_2$, then obviously $(\Delta_1,\Delta_2)$ is a splint of $\Delta$;
one may regard this as a trivial splint.
The next result is straighforward, provided
one uses the axioms for root systems developed, for example, in \cite{Hum}:

\begin{Prp}
Suppose $\iota_1:\Delta_1\hookrightarrow\Delta$ and $\iota_2:\Delta_2\hookrightarrow\Delta$ are metric
embeddings.  Then, if it is non-empty, the intersection of the images of $\iota_1$ and $\iota_2$
is a root system metrically embedded in $\Delta$.
\end{Prp}

Next, it is useful to characterize root systems of type $ADE$ in this context:

\begin{Prp}
Let $\Delta$ be a simple root system.  Then $B_2\hookrightarrow \Delta$ if and only if $\Delta$ has type $BCFG$. 
\end{Prp}

\noindent{\sc Proof.}

Notice that $B_2$ has roots $\alpha$ and $\beta$ such that both $\alpha+\beta$ and $2\alpha+\beta$
are also roots, as does every root system of type $BCFG$.  On the other hand, no root sytem of type
$ADE$ has such a pair of roots.

\hspace{\fill}$\square$

\begin{Cor}
Suppose $\Delta_0$ is simple and $\Delta$ has type $ADE$.  If $\iota:\Delta_0\hookrightarrow\Delta$ is an embedding,
then $\Delta_0$ has type $ADE$ and $\iota$ is metric.
\end{Cor}

\noindent{\sc Proof.}

First use the characterization of root systems of type $ADE$ given above; since $\Delta_0$
is simple, this implies that
$\Delta_0$ has type $ADE$.  Next, assuming all roots have norm 2,
if $\alpha$ and $\beta$ are roots, then $\alpha+\beta$ is a root
if and only if $\left<\alpha,\beta\right>=-1$.  
This shows that $\iota$ is metric.

\hspace{\fill}$\square$

One may now proceed with the various cases.    There are two sections for handling the
special cases, divided according to whether the simple root system has type $ADE$ or $BCFG$.
For each of the infinite sequences $A_r$, $B_r$, $C_r$, and $D_r$, there are two different 
parts to classifying the splints.
First one must determine all splints when $r$ is small, and then
one must show that when $r$ is sufficiently large, the number of splints is severely limited.

\subsection{Root sytems of type $ADE$}

$$A_r$$

Since it has only 1 element, it is obvious that that $A_1$ does not splinter.
Given any $r\geq 3$, there are at least two splints of $A_r$, namely
$(rA_1,A_{r-1})$ and $(A_1+A_{r-1},(r-1)A_1)$.  
One may describe the first splint by writing
$$\Delta_1=\{e_i-e_r:0\leq i\leq r-1\}\hbox{ and }\Delta_2=\{e_i-e_j:0\leq i<j\leq r-1\},$$
and the second by
$$\Delta_1=\{e_i-e_j:0\leq i<j\leq r-1\}\cup\{e_0-e_r\}\hbox{ and }\Delta_2=\{e_i-e_r:1\leq i\leq r-1\}.$$
Evidently the two splints are closely related, one being obtained from the other by
``shifting'' $e_0-e_r$ between the two summands.
These splints coincide when $r=2$, yielding the splint $(A_1,2A_1)$ of $A_2$.
The root system $A_3$ has an additional splint as $(3A_1,3A_1)$, letting
$$\Delta_1=\{e_0-e_1,e_1-e_2,e_0-e_3\}\hbox{ and }\Delta_2=\{e_2-e_3,e_0-e_2,e_1-e_3\}.$$
Notice that in each of these subsets, there is a pair of orthogonal roots.  Thus, it seems appropriate
to denote this splint by $(A_1+D_2,A_1+D_2)$.
The root system $A_4$ has two additional splints as $(A_2+2A_1,A_2+2A_1)$ and $(2A_2,4A_1)$.
One may describe the first by letting
$$\Delta_1=\{e_2-e_3,e_3-e_4,e_2-e_4,e_0-e_3,e_1-e_4\}$$
and
$$\Delta_2=\{e_0-e_1,e_1-e_2,e_0-e_2,e_1-e_3,e_0-e_4\},$$
and the second by
$$\Delta_1=\{e_0-e_1,e_1-e_2,e_0-e_2,e_2-e_3,e_3-e_4,e_2-e_4\}$$
and
$$\Delta_2=\{e_1-e_3,e_0-e_3,e_1-e_4,e_0-e_4\}.$$
Again, due to the presence of pairs of orthogonal roots under the usual inner product, these splints
may be denoted $(A_2+D_2,A_2+D_2)$ and $(2A_2,2D_2)$ respectively.

\begin{Lem}
If $\Delta$ is simple and $\Delta\hookrightarrow A_r$, then $\Delta\cong A_s$ for some
$s\leq r$.
\end{Lem}

\noindent{\sc Proof.}

If one writes the highest root of $A_r$ as a linear combination of simple roots,
then every coefficient is equal to 1.

\hspace{\fill}$\square$

\begin{Prp}

(a) The root system $A_2$ has exactly one splint.
(b) The root system $A_3$ has exactly three splints.
(c) The root system $A_4$ has exactly four splints.
(d) The root system $A_5$ has exactly two splints.

\end{Prp}

\noindent{\sc Proof.}

Parts (a) and (b) are evident.  Parts (c) and (d) follow quickly after observing the following:  If
$A_s\hookrightarrow A_r$ and $A_t\hookrightarrow A_r$ are embeddings with disjoint images, then $s+t\leq r$.
Suppose $r=4$.  It is clear that there are only two splints having $A_3$ as one of the components.
If one of the components is $A_2$, then every other component is isomorphic to either $A_2$ or $A_1$.
By simply listing the cases, one can eliminate all but the two described above.  The argument for
$A_5$ is similar.

\hspace{\fill}$\square$

\begin{Lem}
Suppose $r\geq 3$ and $s\geq 2$.  If $A_r$ has a splint where $A_s$ is a component, then
$A_{r-1}$ has a splint having $A_{s-1}$ as a component.
\end{Lem}

\noindent{\sc Proof.}

Suppose $(\Delta_1,\Delta_2)$ is a splint of $A_r$ with having $\iota:A_s\hookrightarrow\Delta_1$
as a component with $s\geq 2$.  Without loss of generality, one may assume that the roots in the image
of $\iota$ have the form $e_i-e_j$, where $i<j$ and $i,j\in\{0,1,2,...,s-1,r\}$.
Consider restricting the splint to the embedding $\kappa:A_{r-1}\hookrightarrow A_r$
where the image of $\kappa$ consists of roots of the form $e_i-e_j$ where $0\leq i<j\leq r-1$.
Since $\kappa$ is metric and all components of $(\Delta_1,\Delta_2)$ are embedded metrically,
this yields a splint of $A_{r-1}$ having $A_{s-1}$ as a component.

\hspace{\fill}$\square$

\begin{Prp}
Suppose $r\geq 6$.  If $(\Delta_1,\Delta_2)$ is a splint of $A_r$ having $A_s$ as a component,
then $s\in\{1,r-1\}$.
\end{Prp}

\noindent{\sc Proof.}

If $r\geq 6$, then one may argue by contradiction, using the preceeding results.

\hspace{\fill}$\square$

Combining the results of this section, one sees that $A_r$ has exactly two splints when $r\geq 5$.

$$D_r$$

There is only one value of $r\geq 4$ for which $D_r$ splinters, namely $r=4$.  In this case there is only one splint
$(2A_2,2A_2)$, and one may present this by writing
$$\Delta_1=\{e_1-e_2,e_2-e_4,e_1-e_4\}\cup\{e_3+e_4,e_1-e_3,e_1+e_4\}$$
and
$$\Delta_2=\{e_3-e_4,e_2+e_4,e_2+e_3\}\cup\{e_2-e_3,e_1+e_3,e_1+e_2\}.$$

\begin{Lem}
The root system $D_5$ does not splinter.
\end{Lem}

\noindent{\sc Proof.}

Suppose $(\Delta_1,\Delta_2)$ is a splint of $D_5$.  Then one has $|\Delta_1|+|\Delta_2|=|D_5|=20$,
so, without loss of generality, one may assume that $|\Delta_1|\geq 10$.  Also, from an earlier proposition,
every component of $\Delta_1$ is embedded metrically.  There are exactly five root systems with rank not exceeding
5 which have these properties, namely $A_4$, $A_1+A_4$, $D_4$, $A_1+D_4$, and $A_5$.  One must rule out all five
cases.  

Suppose $A_5$ is a component.  Then there is a metric embedding $A_4\hookrightarrow A_5$, 
and one may assume without loss of generality
that the roots in the image of $A_4$ are those of the form $e_i-e_j$ such that $1\leq i<j\leq 5$. 
In particular, the simple roots are $e_i-e_{i+1}$ for $1\leq i\leq 4$.  Every root of $D_5$ not lying
in this image of $A_4$ then has the form $e_i+e_j$.  Since $A_5$ is a component, there is a root of the form
$e_i+e_j$ such that
$$\{e_i-e_{i+1}:1\leq i\leq 4\}\cup\{e_i+e_j\}$$
is a system of simple roots. This is impossible, so $A_5$ cannot be a component.

Next, suppose $A_4$ is a component.  Then the embedding of $A_4$ is metric and, as above, one may assume that the image
consists of all roots of the form $e_i-e_j$, where $1\leq i<j\leq 5$.  Every root not lying
in the image then has the form $e_i+e_j$.  Notice that $A_2$ hence does not embed in the remaining
roots and every other simple component of the splint is isomorphic to $A_1$.  Since there are
10 roots of this form and $A_4$ has rank 4, the total of the ranks of all the simple components is then
14.  Since this exceeds 10, this violates the definition of ``splint''.

An argument in the case that $D_4$ is a component similarly leads to a contradiction.

\hspace{\fill}$\square$

\begin{Prp}
If $r\geq 6$, then $D_r$ does not splinter.
\end{Prp}

\noindent{\sc Proof.}

Suppose $(\Delta_1,\Delta_2)$ is a splint of
$D_r$ where $r\geq 6$.  Then every component is embedded metrically.  Moreover, there
is at least one embedding $D_5\hookrightarrow D_r$ whose image contains elements of both
$\Delta_1$ and $\Delta_2$.  Such an embedding is metric and all the components of the splint
are embedded metrically, so such a splint restricts
to a splint of $D_5$, a contradiction.

\hspace{\fill}$\square$

$$E_r$$

\begin{Prp}
None of the root systems $E_6$, $E_7$, nor $E_8$ splinters.
\end{Prp}

\noindent{\sc Proof.}

Since $D_5$ embeds metrically in $E_r$ for all $r\in\{6,7,8\}$ and $E_r$ has type $ADE$,
one may use the same argument which was used for the cases of $D_r$ with $r\geq 6$.

\hspace{\fill}$\square$

\subsection{Root sytems of type $BCFG$}

$$G_2$$

By inspection, one may quickly establish:

\begin{Prp}
The root system $G_2$ has exactly two splints.
\end{Prp}

These are $(A_2,A_2)$ and $(2A_1,B_2)$.  One obtains the former of these by partitioning $G_2$ into
its long and short roots; thus,
$$\Delta_1=\{2e_1-e_2-e_3,e_1-2e_2+e_3,e_1+e_2-2e_3\}\hbox{ and }\Delta_2=\{e_1-e_2,e_1-e_3,e_2-e_3\}.$$
In this splint, both components are embedded metrically.  One obtains the
latter by considering two of the long roots as copies of $A_1$.  Thus,
$$\Delta_1=\{e_1-2e_2+e_3,e_1+e_2-2e_3\}\hbox{ and }\Delta_2=\{2e_1-e_2-e_3,e_1-e_2,e_1-e_3,e_2-e_3\}.$$

$$F_4$$

Similar to the case of $G_2$, one obtains the totally metric splint $(D_4,D_4)$
when one partitions the roots of $F_4$ into long and short roots.

\begin{Prp}
The root system $F_4$ has exactly one splint.
\end{Prp}

\noindent{\sc Proof.}

One must show that there are no other splints other than that described above.
Suppose $(\Delta_1,\Delta_2)$ is a splint of $F_4$ and consider
the possible cardinalities
of the components $\Delta_1$ and $\Delta_2$.  
One may quickly tabulate the cardinalities of all the root systems, simple or not,
with rank not exceeding 4.  One first finds that $D_4$ is the only such
root system with cardinality 12.   This shows that the only way to splinter $D_4$ into two root systems
each with cardinality 12 is to splinter it as $(D_4,D_4)$.  

At least one of these splints has already been described, so one must now show that
there aren't any others.  A tedious check reveals the existence
of at least one non-metric embedding of $D_4\hookrightarrow F_4$.  This has
simple roots
$$\Sigma=\{e_1-e_4,e_2-e_4,e_3-e_4,2e_4\},$$
and the image consists of 9 short roots 
$$\{e_1\pm e_4,e_2\pm e_4,e_3\pm e_4,e_1+e_2,e_1+e_3,e_2+e_3\}$$
and 3 long roots
$$\{2e_4,e_1+e_2+e_3\pm e_4\}.$$
Moreover, again by some tedious computations, one may see that every other embedding $D_4\hookrightarrow F_4$
which is not one of the two metric embeddings described earlier
has precisely 3 long roots and 9 short roots.
Thus, the remaining set of roots not lying in the image of this embedding, 
which contains 3 short roots and 9 long roots, cannot be
the image of such an embedding.  Therefore, there is only one splint of $F_4$ where
the cardinality of each component is equal to 12.

Next, there are only two root systems having rank not exceeding 4 with
cardinality between 12 and 24, namely $B_4$ and $C_4$, and each of these has
cardinality 16.  Without loss of generality,
assume $\Delta_1$ has cardinality 16.  Then $\Delta_2$ must have cardinality 8, and another 
inspection of the table
reveals that $\Delta_2$ must be isomorphic to $2B_2$.  Thus, the aim now is to rule out the cases that
$(\Delta_1,\Delta_2)$ may be isomorphic to either 
$(B_4,2B_2)$ or $(C_4,2B_2)$.  In order to rule these out, one observes that an embedding of either $B_4$
or $C_4$ in $F_4$ must be metric.  In the former case, $B_4$ contains all the long roots of $F_4$ and
in the latter, $C_4$ contains all the short roots.  In either case, the balance consists
of roots of only one length, and thus cannot contain a copy of $B_2$.

\hspace{\fill}$\square$

$$B_r$$

For any $r\geq 2$, the root system $B_r$ has a splint $(D_r,rA_1)$ into long and short roots, given as
$$\Delta_1=\{e_i\pm e_j:1\leq i<j\leq r\},\hbox{ and }\Delta_2=\{e_i:1\leq i\leq r\}.$$
This splint is distinguished because both $\Delta_1$ and $\Delta_2$ are embedded metrically.
The root system $B_2$ has two additional splints.  The first, $(A_1,A_2)$, may be given by
$$\Delta_1=\{e_1-e_2\}\hbox{ and }\Delta_2=\{e_1,e_2,e_1+e_2\}, $$
and the second, $(2A_1,2A_1)$, may be given by
$$\Delta_1=\{e_1,e_1-e_2\}\hbox{ and }\Delta_2=\{e_2,e_1+e_2\}.$$
The root system $B_3$ has an additional splint as $(3A_1,A_3)$ given by
$$\Delta_1=\{e_1+e_2,e_1+e_3,e_2+e_3\}\hbox{ and }\Delta_2=\{e_1-e_2,e_1-e_3,e_2-e_3,e_1,e_2,e_3\}.$$
There are no other splints of $B_r$ besides these.
The proof depends mainly on the analysis of the $D_r$ root systems and a detailed description
of non-metric embeddings $A_s\hookrightarrow B_r$.  

\begin{Lem}
The root system $C_3$ does not embed in $B_r$ for any $r\geq 2$.
\end{Lem}

\noindent{\sc Proof.}

If $\alpha$, $\beta$, and $\half(\alpha+\beta)$ are roots in $B_r$, then
$\alpha$ and $\beta$ are long and orthogonal to each other, while $\half(\alpha+\beta)$
is short.  Next, $C_3$ contains three roots $\alpha$, $\beta$, and $\gamma$ with the property that
the vectors $\half(\alpha+\beta)$, $\half(\alpha+\gamma)$, and $\half(\beta+\gamma)$ are also roots.
Suppose such roots lie in $B_r$.  All of the short roots in $B_r$ are orthogonal to each other, and so
one has, for instance,
$$\begin{array}{rcl}
0 &=& \left<\alpha+\beta,\alpha+\gamma\right> \\
&=& \left<\alpha,\alpha\right>+
\left<\alpha,\gamma\right>+\left<\beta,\alpha\right>+\left<\beta,\gamma\right> \\
&=& 2+0+0+0 \\
&=& 2, \\
\end{array}$$
a contradiction.

\hspace{\fill}$\square$

\begin{Prp}
If $s\geq 4$, then every embedding $D_s\hookrightarrow B_r$ is metric.
\end{Prp}

\noindent{\sc Proof.}

This may be considered as a corollary to the analysis of the $F_4$ root system.
Every non-metric embedding $D_4\hookrightarrow F_4$
contains exactly 3 long roots and 9 short roots.  Also, there is a metric
embedding $B_4\hookrightarrow F_4$.
However, the root system $B_4$ has only 4 short roots, so there
is no non-metric embedding $D_4\hookrightarrow B_4$.  The cases when $s\geq 5$ follow by restricting
to this case.  

\hspace{\fill}$\square$

\begin{Lem}
If a non-surjective embedding $\iota:\Delta\hookrightarrow B_r$ contains all the long roots of $B_r$, then
$\Delta\cong D_r$.
\end{Lem}

\noindent{\sc Proof.}

This is obvious if $r=2$.  If $r\geq 3$, then $\Delta$ is simple since it has rank $r$ and contains
all the roots of $D_r$.  Since the image of $\iota$ contains all the $r(r-1)$ long roots of $D_r$
and not all of the short roots of $B_r$, there is an integer $k\in\{0,1,2,...,r-1\}$ such that the cardinality
of $\Delta$ is $r(r-1)+k$.  After an examination of the cardinalities of all the simple root systems, one concludes
that $k=0$.

\hspace{\fill}$\square$

Despite the fact that an embedding $\iota:A_s\hookrightarrow B_r$ may not be metric,
describing the image is straightforward:

\begin{Prp}
Suppose $\iota:A_s\hookrightarrow B_r$ is an embedding which is not metric.  Then
the image $\iota(A_s)$ has a splint $(sA_1,A_{s-1})$ where $sA_1$ consists of short
roots of $B_r$ and $A_{s-1}\hookrightarrow B_r$ is a metric embedding.
\end{Prp}

\noindent{\sc Proof.}

If the embedding is not metric, then at least one of the roots is short.  Since there is only one orbit of
roots under the action of the Weyl group for $A_s$, one may assume without loss of generality that
one of the simple roots on one of the ends of its Dynkin diagram is short.  
In this case, one obtains exactly $s$ short roots.  These roots comprise an orthogonal set
and the balance is embedded metrically.

\hspace{\fill}$\square$

\begin{Prp}
The root system $B_4$ has exactly one splint.
\end{Prp}

\noindent{\sc Proof.}

Suppose $(\Delta_1,\Delta_2)$ is a splint of $B_4$.
Since $|B_4|=16$, one may assume without loss of generality that $8\leq|\Delta_1|<16$.
Examining all the root systems with rank less than or equal to 4 and with these cardinalities,
one sees that $\Delta_1$ must be isomorphic to one of $B_3$, $A_1+B_3$, $C_3$, $A_1+C_3$, $2B_2$, or $D_4$.
One must rule out all but one of these possibilities.

It is possible that $\Delta_1$ is isomorphic to $D_4$.  However, in this case, the embedding $D_4\hookrightarrow B_4$
is metric and therefore is the splint $(D_4,4A_1)$ into long and short roots.
Next, $\Delta_1$ cannot be $C_3$ or $A_1+C_3$ because $C_3$ does not embed in $B_r$ for any value of $r$.

Suppose one of the simple components of the splint is $B_3$.  Then the embedding is metric and one
may assume that the image is
$$\{e_1,e_2,e_3,e_1\pm e_2,e_1\pm e_3,e_2\pm e_3\}.$$
One can check that $A_2$ does not embed in the set of 7 roots not lying in this image, and hence
that $B_3$ cannot be a component.

Finally, suppose $\Delta_1\cong 2B_2$.  Then $|\Delta_1|=8$ and hence $|\Delta_2|=8$.  However, there
is only one root system with rank less than or equal to 4 and with cardinality 8, namely $2B_2$.
Since $\Delta_1\cong 2B_2$, each component is metric and hence $\Delta_1$ contains all 4 short
roots of $B_4$.  It is not possible that $\Delta_2$ also contain all 4 short roots of $B_4$,
so $\Delta_1$ cannot be isomorphic to $2B_2$.

\hspace{\fill}$\square$

\begin{Prp}
If $r\geq 5$, then $B_r$ has exactly one splint.
\end{Prp}

\noindent{\sc Proof.}

Suppose $r\geq 5$ and  $(\Delta_1,\Delta_2)$ is a splint of $B_r$
other than $(D_r,rA_1)$.  
Consider restricting this splint to the long roots, that is, to $D_r\hookrightarrow B_r$.
Notice that both $\Delta_1$ and $\Delta_2$ must each contain at least one long root of $B_r$,
so the restriction is not trivial.  
If $\Delta$ is any simple component of the splint
which is not embedded metrically, then $\Delta=A_s$ for some $s$.
However, in this case, the restriction to the long roots is a metric embedding.
Thus, restriction yields a splint of $D_r$ with $r\geq 5$, a contradiction.

\hspace{\fill}$\square$

$$C_r$$

If $r\geq 3$, then there is a splint $(rA_1,D_r)$ given as
$$\Delta_1=\{2e_i:1\leq i\leq r\},\hbox{ and }\Delta_2=\{e_i\pm e_j:1\leq i<j\leq r\}.$$
This splint is distinguished because both components are embedded metrically.
The root system $C_3$ has an additional splint $(A_1+B_2,A_1+A_2)$, which may be presented as
$$\Delta_1=\{e_2-e_3\}\cup\{2e_1,2e_2,e_1\pm e_2\},\hbox{ and }\Delta_2=\{e_2+e_3\}\cup\{2e_3,e_1\pm e_3\}.$$
Evidently, the general conclusion for the case of $C_r$ is similar to the case of $B_r$.
Despite the similarity in the conclusions, however, the proofs differ considerably.
The first part of the proof is to specify which simple components may appear in any splint of $C_r$.
Obviously $A_s$, $C_s$, and $D_s$ may imbed in $C_r$ for certain values of $s$, whereas if $r\geq 4$,
then none of the exceptional root systems may imbed in $C_r$.  One potentiality remains:

\begin{Prp}
If $s\geq 3$, then $B_s$ does not imbed in $C_r$.
\end{Prp}

\noindent{\sc Proof.}

It suffices to prove the case for $s=3$, since $B_3$ embeds in $B_s$ for all $s\geq 3$.
First, if $\alpha$, $\beta$, and $\half(\alpha+\beta)$ are roots in $C_r$, then
$\alpha$ and $\beta$ are long and orthogonal to each other, while $\half(\alpha+\beta)$
is short. Next, $B_3$ contains distinct roots $\{\alpha_1,\alpha_2,\beta_1,\beta_2\}$ 
with the property that the vector
$$\frac{\alpha_1+\beta_1}{2}=\frac{\alpha_2+\beta_2}{2}$$
is also a root.  Suppose such roots lie in $C_r$.  Then one immediately has
$$\left<\alpha_1,\beta_1\right>=\left<\alpha_2,\beta_2\right>=0.$$
One can show that the inner products among the other roots are also zero as follows:
Consider the value $\left<\alpha_1,\alpha_2\right>\in\{\pm 4,0\}$.  One knows
$\left<\alpha_1,\alpha_2\right>\neq 4$ because these roots are not equal.  One knows
$\left<\alpha_1,\alpha_2\right>\neq -4$ because these roots are both positive.  Thus,
$\left<\alpha_1,\alpha_2\right>=0$.  A similar argument holds for the inner products
$\left<\alpha_1,\beta_2\right>$, $\left<\alpha_2,\beta_1\right>$,
and $\left<\beta_1,\beta_2\right>$.  However, this implies that
the ``central'' root $\half(\alpha_1+\beta_2)$ has norm zero, a contradiction.

\hspace{\fill}$\square$

\begin{Prp}
If $s\geq 4$, then every embedding $D_s\hookrightarrow C_r$ is metric.
\end{Prp}

\noindent{\sc Proof.}

Consider the analysis of the $F_4$ root system.  Recall that one has a metric embedding
$C_4\hookrightarrow F_4$.  However, every non-metric embedding $D_4\hookrightarrow F_4$
contains 3 long roots, two having the form $e_1\pm e_i+e_j+e_k$, and the other having the form 
$2e_i$.  Evidently this does not restrict to a non-metric embedding $D_4\hookrightarrow C_4$,
so this establishes the case $s=4$.
If $s>4$, then every embedding $D_s\hookrightarrow C_r$
must restrict to a metric embedding $D_4\hookrightarrow C_r$.  Hence the embedding of $D_s$ is
metric.

\hspace{\fill}$\square$

\begin{Lem}
If $r\geq 4$ and $(\Delta_1,\Delta_2)$ is a splint of $C_r$ where 
every simple component is embedded metrically, then
it is the splint $(rA_1,D_r)$.
\end{Lem}

\noindent{\sc Proof.}

This is routine to establish if $r=4$.  Next, suppose $r\geq 5$.  Let $\iota:D_r\hookrightarrow C_r$
be the unique metric embedding, and restrict the splint to the image of $\iota$.
It was shown earlier that $D_r$ does not splinter when $r\geq 5$, so one may assume
without loss of generality that $\Delta_2$ contains all the short roots of $C_r$.
Since $\Delta_2$ contains all the short roots of $C_r$ while it is not equal to $C_r$, 
it must be simple and therefore isomorphic to $D_r$.

\hspace{\fill}$\square$

It is important to understand the contrapositive of this lemma:  If $(\Delta_1,\Delta_2)$ is a
splint of $C_r$ for which there exists $s\geq 2$ such that $A_s$ is a metrically
embedded component in either $\Delta_1$ or $\Delta_2$, then at least one simple component of either
$\Delta_1$ or $\Delta_2$ is embedded non-metrically.  Moreover, under such circumstances,
such a component is isomorphic to $A_t$, embedded non-metrically, for some $t\geq 2$.

The rest of the proof relies on techniques similar to those used for the cases
of $A_r$, $D_r$, and $B_r$.  However, this case is different because of the potential occurrence
of a non-metrically embedded copy of $A_s$.  Fortunately, it is easy to describe
a non-metric embedding $A_s\hookrightarrow C_r$:

\begin{Lem}
Suppose $\iota:A_s\hookrightarrow C_r$ is a non-metric embedding with $s\geq 2$.
Then 
(a) there is a unique index $i$ such that $2e_i$ lies in the image of $\iota$, and
(b) there is a subset $J\subset\{1,2,3,...,r\}$ having the property that $e_i+e_j$ 
is in the image of $\iota$ if and only if $j\in J$.
\end{Lem}

\noindent{\sc Proof.}

First, if such an embedding contains no long roots, then it restricts to an embedding
$A_s\hookrightarrow D_r$.  Every such embedding is necessarily metric, so the non-metric
embedding $A_s\hookrightarrow C_r$ contains at least one long root, say $2e_i$.
Notice that the root system $A_r$ contains only one orbit of roots under the action of
the corresponding Weyl group, so,
without loss of generality, one may assume that $2e_i$ is the simple root corresponding
to one of the ends of the Dynkin diagram of $A_s$.
If $2e_i$ is a simple root, then all the other simple roots must be short.
This is due to the fact that the coefficients must sum to either 0 or 2.
Were $2e_j$ another simple root, there would be a root in $C_r$ whose coordinate sum
totaled 4.  One can check that the remaining roots must also be short.
Part (b) follows immediately.

\hspace{\fill}$\square$

It is convenient to refer to the set $J$ as the ``distinguished indices''
corresponding to the embedding $\iota:A_s\hookrightarrow A_r$.  The corollary follows
immediately:

\begin{Cor}
Suppose $\iota:A_s\hookrightarrow C_r$ is a non-metric embedding with $s\geq 2$,
$2e_i$ lies in the image of $\iota$, and $J\subset\{1,2,3,...,r\}$ is the set of distinguished indices.
If $j,k\in J$ with $j<k$ and neither $j$ nor $k$ is equal to $i$, then
either $e_j-e_k$ or $e_j+e_k$ (exclusively) is in the image of $\iota$.
\end{Cor}

\begin{Lem}
Suppose $\Delta$ is a simple root system and that
$\iota:C_s\hookrightarrow C_r$ and $\kappa:\Delta\hookrightarrow C_r$ are embeddings.
Then the intersection $\iota(C_s)\cap\kappa(\Delta)$ is a root system.
\end{Lem}

\noindent{\sc Proof.}

First, $\iota$ is necessarily a metric embedding.  Thus,
if $\kappa$ is metric, then the result follows immediately.  The only non-trivial
case occurs if $\Delta$ is $A_t$ for some $t$ and $\kappa$ is non-metric.  Thus, suppose
$\Delta=A_t$ for some $t$ and $\kappa$ is non-metric.  
The form of the roots in the image of $\kappa$ was described in an earlier proposition.
Thus, assume $2e_i$ is the unique long root appearing in the image of $\kappa$ and $I\subset\{1,2,3,...,r\}$
is the set of distinguished indices.
For each element $j\in\{1,2,3,...,r\}$, the restriction to the orthogonal complement of $2e_j$
yields a restriction of $A_s$ to $C_{r-1}$.  
Consider the following three cases:  First, $(2e_i)^\perp\cap \kappa(A_t)$ is a metric
embedding of $A_{t-2}$ in $C_{r-1}$.  Second, if $j\in I$ and $j\neq i$, then 
$(2e_j)^\perp\cap \kappa(A_t)$ is a non-metric embedding of $A_{t-1}$.
Third, if $j\notin I$, then $(2e_j)^\perp\cap \kappa(A_t)=\kappa(A_t)$,
for none of the indices in $\kappa(A_t)$ are equal to $j$.

For the general case, notice that there is a subset $J\subset\{1,2,3,...,r\}$ such that
$2e_i$ lies in the image of $\iota$ if and only if $i\in J$.  One may then apply the above
arguments inductively, starting with $t=r-1$ and proceeding down until $t=|J|$.

\hspace{\fill}$\square$

\begin{Cor}
Suppose $(\Delta_1,\Delta_2)$ is a splint of $C_r$ and $\iota:C_s\hookrightarrow C_r$
is an embedding.  Then the restriction of $(\Delta_1,\Delta_2)$ to the image of $\iota$ is
a splint of $C_s$.
\end{Cor}

\begin{Prp}
If $r\geq 4$, then $C_r$ has exactly one splint.
\end{Prp}

\noindent{\sc Proof.}

Suppose $(\Delta_1,\Delta_2)$ is a splint of $C_r$ which is not equivalent to
$(rA_1,D_r)$.  Then at least one simple component $\Delta$ is embedded non-metrically.  This component
$\Delta$ cannot be $D_s$ or $C_s$ for each of these must embed metrically.  Nor can $\Delta$
be $B_s$ or any exceptional root system for none of these embed in $C_r$.  
The only remaining possibility is that $\Delta\cong A_s$ for some $s\geq 2$.
This immediately leads to a contradiction because then this splint restricts to a
splint of $C_4$ having $A_2$ as a component.

\hspace{\fill}$\square$

\section{Conclusion}

This classification has not addressed the issue of equivalence.  
Indeed, one may define:  
If $\Delta$ is a simple root system with Weyl group $W$, then the splints $(\Delta_1,\Delta_2)$
and $(\Delta_1',\Delta_2')$ of $\Delta$ are ``equivalent'' if there exists $\sigma\in W$ such that
$$\sigma\cdot(\Delta_1\cup(-\Delta_1),\Delta_2\cup(-\Delta_2))=(\Delta_1'\cup(-\Delta_1'),\Delta_2'\cup(-\Delta_2')).$$ 
The general problem is to detect whether or not two splints are equivalent.  As the reader must have noticed by
now, however, splints of simple root systems are scarce and easy to describe, so it would not take much effort to determine
that the table is complete up to equivalence.  The splints of
$A_4$ and $D_4$ probably afford the most non-trivial (and most interesting) cases.

There may be a more efficient way to achieve the classification of splints of simple root systems than that given here, although
this author has not found one. 

\section{Appendix}

Each row in the following table gives a splint $(\Delta_1,\Delta_2)$ of the simple root system $\Delta$.

$$\begin{array}{c c||c|c}
\hbox{type}
& \hspace{0.25 in}\Delta\hspace{0.25 in} & \hspace{0.25 in}\Delta_1\hspace{0.25 in} & \hspace{0.25 in}\Delta _2\hspace{0.25 in} \\
\hline
\hline
\hbox{(i)}     & G_2 & A_2 & A_2 \\
& F_4 & D_4 & D_4 \\
\hline
\hbox{(ii)} & B_r (r\geq 2) & D_r & rA_1 \\
    & C_r (r\geq 3) & rA_1 & D_r \\
\hline
\hbox{(iii)}    & A_r (r\geq 2) & A_{r-1} & r A_1 \\
    & B_2 & A_1 & A_2 \\
\hline
\hbox{(iv)}    & B_2 & 2A_1 & 2 A_1 \\
    & A_3 & A_1+D_2 & A_1+D_2 \\
    & A_4 & A_2+D_2 & A_2+D_2 \\
    & D_4 & 2A_2 & 2A_2 \\
\hline
\hbox{(v)}    & A_4 & 2A_2 & 2D_2 \\
    & A_r (r\geq 3) & (r-1)A_1 & A_1+A_{r-1} \\
\hline
\hbox{(vi)}    & B_3 & 3A_1 & A_3 \\
    & C_3 & A_1+B_2 & A_1+A_2 \\
    & G_2 & 2 A_1 & B_2 \\

\end{array}$$

This table is organized into several types:  
(i) Both $\Delta_1$ and $\Delta_2$ are embedded metrically, and $\Delta_1\cong\Delta_2$. 
(ii) Both $\Delta_1$ and $\Delta_2$ are embedded metrically, but $\Delta_1$ and $\Delta_2$ are not isomorphic.
(iii) Only $\Delta_1$ is embedded metrically; 
(iv) Every simple component is embedded metrically and $\Delta_1\cong\Delta_2$. 
(v) Every simple component is embedded metrically, but $\Delta_1$ and $\Delta_2$ are not isomorphic. 
(vi) At least one simple component of $\Delta_2$ is not embedded metrically.
The notation $D_2$ is used in a few instances to designate a metric embedding of $2A_1$.
The types (i), (ii), and (iii) share the property that $\Delta_1$ corresponds to a Lie subalgebra
of the Lie algebra with root system $\Delta$.

\end{document}